\documentclass[12pt]{amsart}
\usepackage[margin=1in]{geometry}
\usepackage[T1]{fontenc}
\usepackage{lmodern}
\usepackage{amsmath,amssymb,amsthm,mathtools,microtype}
\usepackage{booktabs,array,enumitem}
\usepackage[hidelinks,hypertexnames=false]{hyperref}

\newtheorem{theorem}{Theorem}[section]
\newtheorem{proposition}[theorem]{Proposition}
\newtheorem{lemma}[theorem]{Lemma}

\theoremstyle{definition}

\theoremstyle{remark}

\newcommand{\R}{\mathbb R}
\newcommand{\C}{\mathbb C}
\newcommand{\HH}{\mathbb H}
\newcommand{\OO}{\mathbb O}
\newcommand{\Sym}{\operatorname{Sym}}
\newcommand{\Herm}{\operatorname{Herm}}
\newcommand{\tr}{\operatorname{tr}}
\newcommand{\Tr}{\operatorname{Tr}}
\newcommand{\ip}[2]{\langle #1,#2\rangle}
\newcommand{\lmax}{\lambda_{\max}}
\newcommand{\spnorm}[1]{\lVert #1\rVert}

\title[Singular solutions from Jordan algebras]{Singular solutions from simple Euclidean Jordan algebras}
\author[B.Z. Chu]{BaoZhi Chu}
\address{Department of Mathematics, Northwestern University,
2033 Sheridan Road, Evanston, IL 60208}

\begin{document}

\begin{abstract}
For every finite-dimensional simple Euclidean Jordan algebra $V$ of rank at least three, we realize the function $\operatorname{tr}(x^3)/|x|$ on its traceless subspace $V_0$ as a $C^{1,1}\setminus C^2$ viscosity solution of a uniformly elliptic equation $F(D^2 u)=0$ in $\mathbb{R}^{\operatorname{dim}V_0}$. This unifies and extends earlier constructions from Cartan isoparametric cubics, Hsiang minimal cubics, and traceless real symmetric matrices.
\end{abstract}
\maketitle

\section{Introduction}
We denote by $\Sym(n)$ the space of real
symmetric $n\times n$ matrices. 
Let $F:\Sym(n)\to \R$ be a \emph{uniformly elliptic} operator. Namely, there exist $\Lambda\ge \lambda>0$ such that
\begin{equation}\label{eq:ellipticity}
\lambda\tr N\leq F(M+N)-F(M)\leq\Lambda\tr N\quad
\text{for every } M\in\Sym(n)\text{ and every }N\geq 0. 
\end{equation}

Consider the equation
\begin{equation}\label{eq:PDE}
 F(D^2 u)=0.
\end{equation}

A fundamental question in the study of elliptic equations is whether a continuous viscosity solution $u$ of \eqref{eq:PDE} must be $C^2$. 
The answer in dimension two is positive due to a classical result of Nirenberg \cite{Nirenberg} in the 1950s. Over the past two decades, after a series of works progressively bringing down the dimension,
Nadirashvili--Tkachev--Vl\u{a}du\c{t} \cite{NTV} gave a negative answer to this question in dimensions five and higher. 
The question remains open in dimensions three and four. An important point in this question is that we do not assume convexity or concavity on $F$, since otherwise $u$ must be $C^{2,\alpha}_{\mathrm{loc}}$ by the Evans--Krylov theorem \cite{Evans-82,Krylov-82}.
Some extensions of the Evans--Krylov theorem 
were established by Cabr\'e--Caffarelli \cite{CC-jmpa} and
Caffarelli--Yuan \cite{CY}.
Recently, an extension of Nirenberg's theorem in dimensions three and higher under additional geometric hypotheses was established by the author \cite{Chu}.
For background and related surveys on this question, see, for example, \cite{Caffarelli1995FullyNE, XX-book, Mooney-review, NTVBook, Teixeira-survey} and references therein.

Jordan algebras were introduced by Pascual Jordan in the early
1930s in an attempt to find a new algebraic formulation of quantum
mechanics and later found applications in Lie theory,
differential geometry, and other areas. See \cite{FK, McCrimmon, Koecher} for standard references on Jordan algebras. 
  
In the monograph \cite{NTVBook}, Nadirashvili--Tkachev--Vl\u{a}du\c{t} provide a Jordan algebra approach to extend the solution $P_5(x)/|x|$ in \cite{NTV} to solutions $P_n(x)/|x|$ in dimensions $n=5,8,14$, and $26$. Here $P_n(x)$ is the $n$-dimensional Cartan isoparametric cubic.
These four constructions are in one-to-one correspondence with the isomorphism classes of finite-dimensional simple Euclidean Jordan algebras of rank three.
In this paper, we construct $C^{1,1}\setminus C^2$ viscosity solutions of \eqref{eq:PDE} from the traceless subspace of every finite-dimensional simple Euclidean Jordan algebra of rank at least three.

Let $V$ be a finite-dimensional simple Euclidean Jordan algebra of rank
$r\geq3$. We use $\circ$ to denote its Jordan product. Throughout the paper,
we endow $V$ with the canonical inner product
\begin{equation}\label{equ-canonical-inner-product}
\langle x, y\rangle\coloneqq \tr(x\circ y),\quad x,y \in V,
\end{equation}
and its induced norm
$|x|:=\sqrt{\ip{x}{x}}$, where $\tr$ is the Jordan trace. The terminology and standard facts about Jordan algebras are recalled in \S\ref{sec:background}. 

Consider the traceless subspace of $V$ defined by
\begin{equation}\label{equ-traceless-space}
 V_0:=\{x\in V\mid \tr x=0\}.
\end{equation}
It is known that 
$$
n_0\coloneqq \dim V_0 =(r-1)(1+dr/2),
$$
where the positive integer $d$ is the common
dimension of
off-diagonal subspaces in the Peirce decomposition of $V$.

Let
$$
 U(x):= \tr(x^3)/|x|,
 \quad x\in V_0,
$$
where $U(0)=0$, and $x^3:=(x\circ x)\circ x$. Our main theorem is as follows.

\begin{theorem}\label{thm:main}
Let $V$ be a finite-dimensional simple Euclidean Jordan algebra of rank
$r\geq 3$, and let $\iota:\R^{n_0}\to V_0$ be a linear isometry. 
Then there exists $F:\Sym(n_0)\to\R$ satisfying \eqref{eq:ellipticity} with
$\lambda=1$ and $\Lambda=1+(r+3)(n_0+3)$,
such that $u:=U\circ \iota$ is a viscosity solution of \eqref{eq:PDE} in $\R^{n_0}$. Moreover, $u\in C^{1,1}(\R^{n_0})$ is not twice differentiable at $0$.
\end{theorem}

The fundamental classification result of Jordan--von Neumann--Wigner \cite{JNW} (see also \cite[Theorem V.3.7]{FK}) shows that every simple Euclidean Jordan
algebra $V$ of rank $r\ge 3$ is isomorphic to a space of Hermitian matrices with entries in one of the algebras $\R$, $\C$, $\HH$, or $\OO$.
Here $\C$, $\HH$, and $\OO$ denote the algebras of complex
numbers, quaternions, and octonions, respectively; their real
dimensions are $2$, $4$, and $8$.
All simple Euclidean Jordan algebras of rank $r\ge 3$ are listed as follows.
\begin{center}
\renewcommand{\arraystretch}{1.2}
\begin{tabular}{@{}lllll@{}}
\toprule
$V$ & $r$ & $d$ & $n_0$ & First several $n_0$'s\\
\midrule
$\Sym(r,\R)$ & $r\geq3$ & $1$ & $r(r+1)/2-1$ & $5,9,14,20,27,35,44,54,\ldots$ \\
$\Herm(r,\C)$ & $r\geq3$ & $2$ & $r^2-1$ & $8,15,24,35,48,63,80,99,\ldots$\\
$\Herm(r,\HH)$ & $r\geq3$ & $4$ & $2r^2-r-1$ & $14,27,44,65,90,119,152,189,\ldots$ \\
$\Herm(3,\OO)$ & $3$ & $8$ & $26$ & 26\\
\bottomrule
\end{tabular}
\end{center}

For every possible $n_0$ listed in the above table,
Theorem \ref{thm:main} constructs a $C^{1,1}\setminus C^2$ viscosity solution $u$ of a uniformly elliptic equation \eqref{eq:PDE} in $\R^{n_0}$.

When $r=3$, the functions in Theorem \ref{thm:main} coincide, up to linear isometries and nonzero constant factors, with solutions constructed from Cartan isoparametric cubics by Nadirashvili--Tkachev--Vl\u{a}du\c{t}
 in dimensions $5,8,14,26$; see \cite{NTV} and \cite[Chapter 5]{NTVBook}.

When $V=\Sym(r,\R)$, the conclusion of
Theorem \ref{thm:main} was recently established
by the author \cite{Chu}.
In particular, when $V=\Sym(4,\R)$, the function
in Theorem \ref{thm:main} coincides, up to a linear
isometry and a nonzero constant factor, with the
nine-dimensional determinant-quotient solution
recorded in \cite[Remark 4.1.5]{NTVBook} and \cite[p.1008]{S}.
Its numerator is equivalent to the
nine-dimensional Hsiang minimal cubic; see the discussion
in \cite{Chu} and the references therein.

Beyond these previously established cases,
Theorem \ref{thm:main} extends the realization of $U$ as a singular viscosity solution of a uniformly elliptic
equation to the families
$V=\Herm(r,\C)$ and $V=\Herm(r,\HH)$ for every $r\geq4$.
It also provides a unified approach for all simple Euclidean
Jordan algebras of rank at least three, with explicit ellipticity
bounds.

The rank two case cannot be included in Theorem \ref{thm:main}, since then $U\equiv 0$ does not lead to singular solutions.
For the rank-three cases,
Nadirashvili--Tkachev--Vl\u{a}du\c{t} also realize the Cartan quotients as solutions of uniformly elliptic
Hessian equations, whose operators depend only on
the eigenvalues of their matrix arguments.
 However, this cannot be simply extended to the general rank cases. Indeed, as observed in \cite[Remark 4.1.5]{NTVBook},
when $V=\Sym(4,\R)$, the function
$u$
cannot solve a uniformly elliptic Hessian equation $F(D^2 u)=0$.

The proof of Theorem \ref{thm:main} extends the approach developed
in \cite[Section~5]{Chu} for traceless real symmetric matrices
to simple Euclidean Jordan algebras of arbitrary rank.
Using the spectral theorem and the Peirce decomposition,
we establish Jordan algebra analogues of the matrix inequalities
in \cite{Chu} for traceless real symmetric matrices.
We then
use two test directions
in the traceless part of a rank-two Jordan subalgebra
to
obtain a key estimate controlling
the trace of a Hessian difference by its largest eigenvalue. 
Our approach applies uniformly at arbitrary rank and avoids
the explicit Hessian eigenvalue calculations used in the
rank-three treatment of the Cartan cubic solutions
in \cite[Section 5.1]{NTVBook}.

The paper is organized as follows.
Section \ref{sec:traceless} establishes the algebraic inequalities
and properties of the traceless subspace $V_0$ needed in the proof.
Section \ref{sec:hessian} proves the key estimate for differences
of Hessians.
Section \ref{sec:operator} proves Theorem~\ref{thm:main}.

\section{Preliminary facts for Jordan algebras}\label{sec:background}

In this section, we review standard facts from Jordan algebras; see \cite{FK} for more details.

A \emph{Jordan algebra} is a real vector space $V$ with a bilinear
product $\circ$ satisfying
$$
 x\circ y=y\circ x\quad\text{and}\quad
 x\circ(x^2\circ y)=x^2\circ(x\circ y),\quad \text{for every }x,y \in V.
$$
Here $x^2=x\circ x$ and the second identity is called the Jordan identity. For each $x\in V$, the multiplication by $x$ defines a linear operator $L(x): V\to V$ by $L(x)y\coloneqq x\circ y$. Therefore, the second defining property of a Jordan algebra says that endomorphisms $L(x)$ and $L(x^2)$ commute, i.e. $L(x)L(x^2)=L(x^2)L(x)$.

In general, a Jordan algebra need not be associative, while it is always \emph{power associative}. Namely, for every $x\in V$, write $x^1=x$ and define inductively $x^n=x\circ x^{n-1}$, $n=2,3,\dots$. Then $x^p\circ x^q=x^{p+q}$ for all positive integers $p,q$.

Here and throughout the paper, we assume that a Jordan algebra is finite dimensional (as a vector space) and has the identity element $e$ (i.e. $e\circ x=x$, $\forall x$). For such a Jordan algebra $V$, it is \emph{Euclidean} if there exists an inner product on $V$ denoted by $\langle\cdot,\cdot\rangle$ such that $$\langle L(x)u,v\rangle=\langle u, L(x)v\rangle \quad \text{for every }x,u,v\in V.$$ Equivalently, every $L(x)$ is self-adjoint with respect to $\langle\cdot,\cdot\rangle$.

Given a Jordan algebra $V$, recall that an ideal of $V$ is a linear subspace $I$ with $V\circ I\subset I$. We say $V$ is \emph{simple} if $V\ne \{0\}$ and its only ideals are $\{0\}$ and $V$. 

As a model example, the space of $r\times r$ real symmetric matrices
 $\Sym(r,\R)$ with 
 $x\circ y =\tfrac12(xy+yx)$ defines a simple Euclidean Jordan algebra. Its inner product is given by the standard Frobenius product $\tr(xy)$ and its identity element is the identity matrix $I_r$.

\subsection{The spectral theorem}
Let $V$ be a Euclidean Jordan algebra.
We say $c\in V$ is
an \emph{idempotent} if $c^2=c$. Idempotents $c,d$ are \emph{orthogonal} if
$c\circ d=0$, or equivalently $\langle c,d\rangle=0$. An idempotent is \emph{primitive} if it is nonzero and it cannot be written as the sum of two nonzero orthogonal idempotents. A \emph{Jordan frame} is a
family of orthogonal primitive idempotents whose sum is $e$. The following is the spectral theorem for $V$; see \cite[\S III.1]{FK}.

\begin{theorem}\label{thm:spectral}
    Let $V$ be a Euclidean Jordan algebra. Then for $x\in V$ there exist a Jordan frame $c_1,\dots,c_r$ and real numbers $\lambda_1,\dots,\lambda_r$ such that 
    $$x=\sum_{j=1}^r \lambda_j c_j.$$
    Moreover, numbers $\lambda_j$ with their multiplicity are uniquely determined by $x$. Furthermore, if $\lambda_1,\dots,\lambda_r$ are distinct, then $c_1,\dots,c_r$ are also uniquely determined by $x$. 
\end{theorem}
We call the numbers $\lambda_j$ the \emph{eigenvalues} of $x$, and define the \emph{Jordan trace} of $x$ by
$$\tr x\coloneqq \sum_{j=1}^r \lambda_j.$$
Then $\tr: V\to \R$ is linear. It follows from the spectral theorem that every Jordan frame of $V$ has the same number $r$ of elements, called the \emph{rank} of $V$. The spectral theorem also implies the following useful fact:
\begin{lemma}\label{lem:spectral-facts}
    Let $V$ be a Euclidean Jordan algebra of rank $r$, and let $x\in V$ with eigenvalues $\lambda_1,\dots,\lambda_r$. The following hold.
    \begin{enumerate}
        \item[(i)] $\tr(x^k)=\sum_{j=1}^r \lambda_j^k$ for $k=1,2,\dots$.
        \item[(ii)] Assume $x$ is an idempotent. Then $x$ is primitive if and only if $\tr x=1$. 
    \end{enumerate}
\end{lemma}
 
If $V$ is simple, then there exists some positive constant $C$ such that $\langle x,y\rangle =C \tr(x\circ y)$ for $x,y \in V$; see \cite[Proposition III.4.1]{FK}. Throughout the paper, we adopt the canonical inner product \eqref{equ-canonical-inner-product} on $V$,
and its induced norm $|x|:= \sqrt{\tr(x^2)}$. Observe that every Jordan frame is orthonormal for this inner product, and $|e|^2=r$. 
We also use the spectral norm defined by
$$
 \spnorm{x}:=\max_j|\lambda_j|,
$$
where $\lambda_j$ are the eigenvalues of $x$.

\subsection{The Peirce decomposition}
Let $V$ be a Euclidean Jordan algebra with a fixed Jordan frame $c_1,\ldots,c_r$. For every idempotent $c\in V$, the Jordan identity implies that the only possible eigenvalues of $L(c)$ are $0,1,1/2$. 
This induces the following subspaces of $V$ 
$$
V_{ii}:= \R c_i \quad \text{and}\quad 
 V_{ij}:=\{x\in V \mid c_i\circ x=c_j\circ x=\tfrac12 x\} \quad \text{for }i\ne j.      
$$
 For every $x\in V_{ij}$,
it holds that $c_k\circ x=0$ for each $k\notin \{i,j\}$. 
The following Peirce decomposition allows us to compute quadratic forms
by components, just as one expands matrix entries; see \cite[\S IV.1-2]{FK}.

\begin{theorem}\label{thm:peirce}
Let $V$ be a Euclidean Jordan algebra
with a Jordan frame $c_1,\ldots,c_r$.
Then there is an orthogonal direct sum
$$
 V=\bigoplus_{1\le i\le j\le r}V_{ij}.
$$
Moreover, $V_{ij}\circ V_{k\ell}=\{0\}$ if $\{i,j\}\cap\{k,\ell\}=\varnothing$, 
 $V_{ij}\circ V_{jk}\subset V_{ik}$ if $i\ne k$,  and $V_{ij}\circ V_{ij}\subset V_{ii} + V_{jj}$ for every $i,j$. 
\end{theorem}

If we further assume that $V$ is simple,
then
the off-diagonal spaces $V_{ij}$ with $i\ne j$ 
have a common dimension
$d>0$, independent of the choice of Jordan frame; see \cite[Proposition IV.2.3 and Corollary IV.2.6]{FK}.
Hence, 
$\dim V=r+\frac{d\,r(r-1)}2$.

\section{The traceless space \texorpdfstring{$V_0$}{V0}}\label{sec:traceless}
Let $V$ be a finite-dimensional simple Euclidean Jordan algebra of rank
$r\geq 2$. In this section, we study the algebraic properties of its traceless subspace $V_0$ defined in \eqref{equ-traceless-space}.
Since $\tr x=\ip{x}{e}$ for $x\in V$, we have $$V_0=e^\perp \coloneqq \{x\in V\mid \ip{x}{e}=0\}.$$
In particular, $n_0\coloneqq \dim V_0=\dim V-1$. 
Observe that in a Peirce decomposition, every off-diagonal subspace $V_{ij}$ is contained in $V_0$.

Note that $V_0$ is a linear subspace, but is not a Jordan subalgebra, i.e. $V_0$ is not closed under the Jordan product. Indeed, for any nonzero $x\in V_0$, we have $\tr(x^2)=|x|^2>0$ and thus $x^2\notin V_0$.

\subsection{Some algebraic inequalities}
\label{sec:algebra}

The first lemma is an extension to the Euclidean Jordan algebra setting of some standard facts for real symmetric matrices.

\begin{lemma}\label{lem:spectral-bound}
Let $V$ be a Euclidean Jordan algebra of rank $r\geq 2$. Then the following hold.
    \begin{enumerate}
        \item[(i)] 
$\spnorm{x}^2\leq\frac{r-1}{r}|x|^2$ for every $x\in V_0$.
        \item[(ii)]  $0\leq\ip{c}{x^2}\leq\spnorm{x}^2$ for every primitive idempotent $c$ and every $x\in V$.
    \end{enumerate}
\end{lemma}

\begin{proof}
For (i), write the eigenvalues of $x$ as $\lambda_1,\ldots, \lambda_r$.
By Lemma \ref{lem:spectral-facts}(i), we have 
$|x|^2=\tr(x^2)=\sum_i\lambda_i^2$. Since $x$ is traceless, we have  $\lambda_i=-\sum_{j\ne i}\lambda_j$ for every $i$. Thus the Cauchy--Schwarz inequality gives
$$ \textstyle
 \lambda_i^2\leq(r-1)\sum_{j\ne i}\lambda_j^2
             =(r-1)(|x|^2-\lambda_i^2)\quad \text{for every }i.
$$
This proves (i). 

For (ii), since $c$ is an idempotent, the only possible
eigenvalues of $L(c)$ are $0,1$, and $1/2$. Hence, $\ip{c}{x^2}=\ip{L(c)x}{x}\ge 0$. Now we write $x=\sum_i \lambda_i c_i$ in a Jordan frame. Then 
$$  \textstyle
 \ip{c}{x^2}=\sum_i\lambda_i^2\ip{c}{c_i}
 \leq\spnorm{x}^2\sum_i\ip{c}{c_i}
 =\spnorm{x}^2\ip{c}{e}=\spnorm{x}^2.
$$
The last equality uses the primitivity of $c$ and 
Lemma \ref{lem:spectral-facts}(ii). 
\end{proof}

The following lemma is the Jordan form of the matrix inequalities in
\cite[Lemma 5.1]{Chu}.

\begin{lemma}\label{lem:main-algebra}
Let $V$ be a Euclidean Jordan algebra of rank $r\geq2$, let $x\in V_0$, and let $s$ be the sum of the first and second largest eigenvalues of $x$.
Then, for every $y\in V_0$,
\begin{equation}\label{eq:algebra-general}
 \ip{x}{y^2}\leq s|y|^2.
\end{equation}
If we further assume $\ip{x}{y}=0$, then
\begin{equation}\label{eq:algebra-orthogonal}
 \ip{x}{y^2}\leq\frac{r-1}{r}s|y|^2.
\end{equation}
\end{lemma}

\begin{proof}
Write $x=\sum_{i=1}^r \lambda_i c_i$ in a Jordan frame $c_1,\ldots,c_r$. Without loss of generality, we may assume $\lambda_1\ge \ldots\ge \lambda_r$. Then $s=\lambda_1+\lambda_2$. 

We first prove \eqref{eq:algebra-general}. 
Using $e=\sum_i c_i$ and Lemma \ref{lem:spectral-bound}(ii), we obtain
\begin{align}\label{eq:row-estimate}
 \ip{x}{y^2}
 &=  (\lambda_1-\lambda_2)\ip{c_1}{y^2}
 +\lambda_2|y|^2
 -\sum_{j=2}^r(\lambda_2-\lambda_j)\ip{c_j}{y^2}
 \notag\\
 &\leq
(\lambda_1-\lambda_2)\ip{c_1}{y^2} +  \lambda_2|y|^2.
\end{align}
Applying Lemma \ref{lem:spectral-bound}(i)-(ii) to the first term gives
\begin{equation}\label{eq:theta-estimate}
 \ip{x}{y^2}
 \leq
 \left(s-\frac{\lambda_1+(r-1)\lambda_2}{r}\right)|y|^2.
\end{equation}
Since $x$ is traceless, we obtain
$\lambda_1+(r-1)\lambda_2
 =\sum_{j=2}^r(\lambda_2-\lambda_j)\geq 0$, and therefore 
\eqref{eq:algebra-general} is proved.

Now we assume $\ip{x}{y}=0$ and prove \eqref{eq:algebra-orthogonal}. The proof is divided into two cases.

\emph{Case 1. $\lambda_2\geq0$.}
Since
$\lambda_1+(r-1)\lambda_2
 =s+(r-2)\lambda_2\geq s$, 
the desired inequality follows from \eqref{eq:theta-estimate}.

\emph{Case 2. $\lambda_2<0$.}
Write the Peirce decomposition
$$
 y=\sum_i y_i c_i+\sum_{i<j}y_{ij} \quad \text{for }
 y_i\in \R \text{ and } y_{ij}\in V_{ij}.
$$
In particular,
$|y|^2=\sum_i y_i^2+\sum_{i<j}|y_{ij}|^2$.

From $\tr y=0$ and $\ip{x}{y}=0$, we obtain
$\sum_i y_i=0$ and $\sum_i\lambda_i y_i=0$,
and consequently
\[
 (\lambda_1-\lambda_2)y_1
 =\sum_{j=2}^r(\lambda_2-\lambda_j)y_j.
\]
The coefficients on the right are nonnegative, and
$0\leq\sum_{j=2}^r(\lambda_2-\lambda_j)
 =\lambda_1+(r-1)\lambda_2
 <\lambda_1-\lambda_2$.
Thus Cauchy--Schwarz gives
 \begin{align*}
 (\lambda_1-\lambda_2)^2y_1^2
 \leq
 \left(\sum_{j=2}^r(\lambda_2-\lambda_j)\right)^2
 \sum_{j=2}^r y_j^2   \leq
 (\lambda_1-\lambda_2)^2\sum_{j=2}^r y_j^2.
 \end{align*}
Since $\lambda_1-\lambda_2>0$, it follows that
$y_1^2\leq\sum_{j=2}^r y_j^2$.
The orthogonality of the Peirce decomposition therefore implies
$$
 2\ip{c_1}{y^2} =2\ip{c_1\circ y}{y}
 =2y_1^2+\sum_{j=2}^r|y_{1j}|^2
 \leq \sum_{j=1}^r y_j^2 +\sum_{j=2}^r |y_{1j}|^2 \le |y|^2.
$$
Substituting this into \eqref{eq:row-estimate}, we obtain
$$
 \ip{x}{y^2}
 \leq
\frac{s}{2}|y|^2
 \leq\frac{r-1}{r}s|y|^2.
$$
Here we used $s\ge 0$ which can be seen from \eqref{eq:algebra-general} by taking $y=c_1-c_2$. Hence, \eqref{eq:algebra-orthogonal} is proved.
\end{proof}

\subsection{The traceless subspace of a rank two Jordan subalgebra}
Let $V$ be a Euclidean Jordan algebra with a fixed Jordan frame $c_1,\ldots,c_r$. Consider 
$$
W_{12} \coloneqq \R (c_1-c_2) \oplus V_{12}.
$$

Clearly, $W_{12}\subset V_0$ is a linear subspace, but it is not a Jordan subalgebra since $(c_1-c_2)^2=c_1+c_2$. When $V$ is simple, $\dim W_{12}=1+d\ge 2$.

On the other hand, one can recognize $W_{12}$ as the traceless part of the Jordan subalgebra $$E_1:=V_{11}\oplus V_{22}\oplus V_{12}=\R(c_1+c_2)\oplus W_{12}.$$ Namely, $W_{12}=\{x\in E_1 \mid \tr x=0\}$. Observe that $\{c_1, c_2\}$ is a Jordan frame of $E_1$ and thus $E_1$ is rank two.

\begin{lemma}\label{lem:rank-two}
Let $V$ be a Euclidean Jordan algebra with a Jordan frame $c_1,\ldots,c_r$, and let $w\in W_{12}$ satisfy $|w|=1$.
Then $w^2=2^{-1}(c_1+c_2)$ and 
$$
 2\sqrt2\,|\ip{v^2}{w}|\leq |v|^2- \ip{v}{u}^2,
$$
for every $v\in w^{\perp}$ and every unit $u\in E_1$. 
\end{lemma}
Here $w^\perp\coloneqq \{z\in V \mid \ip{z}{w}=0\}$ is the orthogonal complement of $w$ in $V$.

\begin{proof}
Denote $c_{1,2}\coloneqq c_1+c_2$. 
Write
$w=t(c_1-c_2)+z$ for some $t\in \R$ and $z\in V_{12}$.
Since $(c_1-c_2)\circ z=0$, we have $w^2=t^2c_{1,2}+z^2$.  Since $z\in V_{12}$, we have $z^2\in V_{11}\oplus V_{22}$ and  
$\ip{z^2}{c_i}=\ip{z\circ c_i}{z}=2^{-1}|z|^2$ for $i=1,2$. It follows that
$$z^2 =2^{-1}|z|^2 c_{1,2}.$$
By the orthogonality of the Peirce decomposition, we have $|w|^2 =2t^2 +|z|^2$. Hence, 
\begin{equation}\label{equ-u-square}
    w^2 = 2^{-1}|w|^2 c_{1,2}. 
\end{equation}
This identity holds for all $w\in W_{12}$.
Since $|w|=1$, it proves the first assertion.

Now we prove the asserted inequality. 
The eigenspaces of $L(c_{1,2})$ corresponding to $1,1/2,0$ are
\begin{align*}
 E_1 := V_{11}\oplus V_{22}\oplus V_{12},\quad 
 E_{1/2}:=\bigoplus_{j=3}^r(V_{1j}\oplus V_{2j}),\quad \text{and}\quad
 E_0:=\bigoplus_{3\leq i \le j\leq r}V_{ij}.
\end{align*}
Write $v=v_1+v_{1/2}+v_0$ corresponding to this orthogonal decomposition. 
By the ``Moreover" part of Theorem \ref{thm:peirce}, the operator $L(w)$ preserves $E_1$ and $E_{1/2}$, and $L(w)$ vanishes on $E_0$.
Therefore, 
\begin{equation}\label{eq:k-split}
 \ip{v^2}{w}=\ip{L(w)v}{v}
            =\ip{v_1^2}{w}+\ip{L(w)v_{1/2}}{v_{1/2}}.
\end{equation}
To compute the first term, write $v_1=\tau c_{1,2}+z$ for some $z\in W_{12}$. Using $c_{1,2}\circ z=z$ and $z^2=2^{-1}|z|^2 c_{1,2}$ (apply \eqref{equ-u-square} to $z$), we obtain
\[
 v_1^2=(\tau^2+|z|^2/2)c_{1,2}+2\tau z.
\]
Since $\ip{c_{1,2}}{w}=0$ and $\ip{v_1}{w}=\ip{v}{w}=0$,
we have
$\ip{z}{w}=0$, and thus
\begin{equation}\label{equ-b1-orthogonal-to-u}
\ip{v_1^2}{w}=0.
\end{equation}

Next we prove that the only possible eigenvalues of $L(w)$ on $E_{1/2}$ are $\pm 1/(2\sqrt{2})$. In fact, we prove $L(w)^2 y =8^{-1}y$ for every $y\in E_{1/2}$. Indeed, fix any such $y$ and consider $x=w+ty$ for $t\in \R$. Then $x^2=2^{-1}c_{1,2}+2t w\circ y +O(t^2)$. Comparing the coefficients of $t$ on both sides of the Jordan identity 
$$x^2\circ (x\circ c_{1,2})= x\circ (x^2 \circ c_{1,2}),$$
we obtain
$8^{-1}y + 2 w\circ (w\circ y)= 4^{-1}y +w \circ (w\circ y)$. Hence, $L(w)^2 y=8^{-1}y$.

From \eqref{eq:k-split}, \eqref{equ-b1-orthogonal-to-u}, and the fact we just proved, we obtain 
\[
 |\ip{v^2}{w}|\leq\frac{1}{2\sqrt2} |v_{1/2}|^2.
\]
Finally, since $u\in E_1$ and $|u|=1$, we have
$\ip{v}{u}^2\leq|v_1|^2$, and thus $|v_{1/2}|^2\leq|v|^2-\ip{v}{u}^2$.
The Lemma is proved.
\end{proof}

\subsection{The orthogonal projection \texorpdfstring{$\pi_0$}{pi0}}\label{sec-3-3}
Let $\pi_0: V\to V_0$ be the orthogonal projection onto $V_0$. 
By $|e|^2=r$, we have
$$
 \pi_0 x=x-\frac{\tr x}{r}e\quad \text{for } x\in V.
$$

For $x\in V$, the operator $\pi_0L(x)|_{V_0}$ is self-adjoint on $V_0$. Indeed,
for $v,w\in V_0$,
\[
 \ip{\pi_0L(x)v}{w}=\ip{L(x)v}{w}
 =\ip{v}{L(x)w}=\ip{v}{\pi_0L(x)w}.
\]
We call $\pi_0L(x)|_{V_0}$ the \emph{compression} of $L(x)$ to $V_0$.

For $E=V$ or $V_0$, we denote by $\Tr_E A$ the ordinary
trace of a linear operator $A:E\to E$. Namely, for an
orthonormal basis $E_1,\ldots,E_{\dim E}$ of $E$,
$$ \textstyle
 \Tr_E A:=\sum_{\alpha=1}^{\dim E}
          \langle AE_\alpha,E_\alpha\rangle.
$$
It is elementary to see that this definition is independent of the choice of orthonormal basis.
We omit the subscript $E$ when the underlying space is clear.

\begin{lemma}\label{lem:operator-traces}
Let $V$ be a simple Euclidean Jordan algebra of rank $r\geq 2$.
Then
$$
 \Tr_{V_0}\bigl(\pi_0L(x)|_{V_0}\bigr)=\frac{n_0}{r}\tr x \quad \text{for every } x\in V. 
$$
\end{lemma}

\begin{proof}
By the spectral theorem in Theorem~\ref{thm:spectral}, we may write $x=\sum_i \lambda_i c_i$ in a
Jordan frame. Then
$$
 L(x)c_i=\lambda_i c_i,\quad \forall i,  \quad \text{and}\quad 
 L(x)z= 2^{-1}{(\lambda_i+\lambda_j)}z,
 \quad\forall z\in V_{ij},\ i\ne j.
$$
By the Peirce decomposition in Theorem~\ref{thm:peirce} and the dimension formulae for $V_{ij}$, we obtain 
$$ 
 \Tr_V L(x)=\sum_i \lambda_i+\frac d2\sum_{i<j}(\lambda_i+\lambda_j)
 =\left(1+\frac{d(r-1)}2\right)\tr x
 =\frac{n_0+1}{r}\tr x.
$$
This formula can be found in
\cite[Proposition III.4.2(i)]{FK}; we include a proof here for the reader's convenience.

Let $E_1,\ldots,E_{n_0}$ be an orthonormal basis of $V_0$.
Since $V_0=e^\perp$ and $|e|^2=r$, adjoining $e/\sqrt r$ gives
an orthonormal basis of $V$. Therefore, we obtain
\begin{align*}
 \Tr_{V_0}\bigl(\pi_0L(x)|_{V_0}\bigr)
 &=\sum_{\alpha=1}^{n_0}\ip{L(x)E_\alpha}{E_\alpha}=\Tr_V L(x)-\ip{L(x)(e/\sqrt r)}{e/\sqrt r}\\
 &=\frac{n_0+1}{r}\tr x-\frac1r\ip{x}{e}
 =\frac{n_0}{r}\tr x.
\end{align*}
The Lemma is therefore proved. 
\end{proof}

\section{Hessian analysis on \texorpdfstring{$V_0$}{V0}}\label{sec:hessian}
In this section, we perform calculus intrinsically on $V_0$.
For a smooth function $f:V_0\to\R$, and for $x\in V_0$, define
$$
 Df(x)[v]:=\left.\frac{d}{dt}f(x+tv)\right|_{t=0}\quad\text{for } v\in V_0, \quad \text{and} 
 $$
 $$
 D^2f(x)[v,w]
 :=\left.\frac{\partial^2}{\partial s\,\partial t}
       f(x+sv+tw)\right|_{s=t=0} \quad \text{for }v,w\in V_0.
$$
Then $Df(x)$ is a linear functional on $V_0$ and $D^2f(x)$ is a symmetric
bilinear form on $V_0$. The latter is represented by a unique self-adjoint
operator $D^2_{V_0}f(x):V_0\to V_0$ defined by 
 $$
 \ip{D^2_{V_0}f(x)v}{w} = D^2f(x)[v,w]
 \quad \text{for }v,w\in V_0.
$$
We write $\Delta_{V_0}f:=\Tr_{V_0}D^2_{V_0}f$ for the Laplacian on $V_0$.

 For $x\in V_0\setminus\{0\}$, let
$$
 H(x):=D^2_{V_0}U(x),
$$
which is a self-adjoint operator on $V_0$.
The following is the main estimate in the proof of Theorem \ref{thm:main}. When $V=\Sym(r,\R)$, it was proved in \cite{Chu}.

\begin{proposition}\label{prop:separation}
Let $V$ be a simple Euclidean Jordan algebra of rank $r\geq 2$.
Then, for every unit $x,y\in V_0$,
\begin{equation}\label{eq:absolute-trace}
 \kappa\lmax(H(x)-H(y))\geq|\Tr_{V_0} (H(x)-H(y))|,
\end{equation}
where $\kappa:=(r+3)(n_0+3)$.
\end{proposition} 

Here $\lmax(A)$ denotes the largest eigenvalue of a self-adjoint
operator $A$ on $V_0$.
For convenience, in the rest of this section, let
$$
 p(x):=\tr(x^3)=\ip{x^2}{x},\quad x\in V_0.
$$
Then $U(x)=p(x)/|x|$ for $x\ne0 $. 
The proof of Proposition \ref{prop:separation} needs the following lemma.

\begin{lemma}\label{lem:H-formula}
Let $V$ be a Euclidean Jordan algebra of rank $r\geq 2$.
For every unit $x,v\in V_0$, 
$$
 \ip{H(x)v}{v}
 =6\ip{x}{v^2}-6\ip{x^2}{v}\ip{x}{v}
   +3p(x)\ip{x}{v}^2-p(x).
$$
If we further assume that $V$ is simple, then
$\Tr_{V_0} H(x)=-(n_0+3)p(x)$.
\end{lemma}

\begin{proof}
Let us first compute derivatives of $p$:
\begin{equation}\label{eq:p-derivatives}
 Dp(x)[v]=3\ip{x^2}{v} \quad \text{and}\quad
 D^2p(x)[v,w]=6\ip{x\circ v}{w} \quad \text{for } x,v,w\in V_0.
\end{equation}

For $|x|=1$ and $v,w\in V_0$,  a direct computation yields
$$
 D(|x|^{-1})[v]=-\ip{x}{v} \quad \text{and} \quad
 D^2(|x|^{-1})[v,w]
 =3\ip{x}{v}\ip{x}{w}-\ip{v}{w}.
$$
Using $U(x)=p(x)|x|^{-1}$, the product rule and
\eqref{eq:p-derivatives} give
\begin{align*}
 \ip{H(x)v}{w}
 ={}&6\ip{x\circ v}{w}
 -3\ip{x^2}{v}\ip{x}{w}
 -3\ip{x^2}{w}\ip{x}{v}\\
 &+3p(x)\ip{x}{v}\ip{x}{w}-p(x)\ip{v}{w}.
\end{align*}
Taking $w=v$, using $|v|=1$ proves the first formula.

To compute the trace, take a sum on the identity we just proved over $v$ equal to the
vectors of an orthonormal basis of $V_0$. By \eqref{eq:p-derivatives}, the first term
gives $\Delta_{V_0}p(x)$. The second term gives
$-6\ip{x^2}{x}=-6p(x)$ by the bilinearity of the inner product. The last two terms give 
$3p(x)|x|^2=3p(x)$ and $-n_0p(x)$, respectively.
Therefore, 
$$
 \Tr_{V_0} H(x)=\Delta_{V_0}p(x)-(n_0+3)p(x).
$$

When $V$ is simple, the second assertion follows from the above and $\Delta_{V_0}p=0$ on $V_0$. To see the harmonicity of $p$, let $\pi_0$ be the orthogonal projection onto $V_0$ as in \S\ref{sec-3-3}. Then the Hessian formula in \eqref{eq:p-derivatives} implies 
$D^2_{V_0}p(x)=6\pi_0L(x)|_{V_0}$.
From this and Lemma \ref{lem:operator-traces}, we obtain for every $x\in V_0$ that
$$
 \Delta_{V_0}p(x)
 =6\Tr_{V_0}\bigl(\pi_0L(x)|_{V_0}\bigr)
 ={6n_0}r^{-1}\tr x=0. 
$$
The Lemma is proved.
\end{proof}

Now we are ready to give the following
\begin{proof}[Proof of Proposition~\ref{prop:separation}]
When $r=2$, the eigenvalues of any traceless element are of the form $\lambda$, $-\lambda$, so 
Lemma \ref{lem:spectral-facts}(i) gives $p\equiv 0$ on $V_0$ and the desired estimate is clear. The following proof treats the case $r\ge 3$.

Fix unit $x,y\in V_0$, and let 
$$T:=H(x)-H(y).$$ 
For the rest of this proof, we write $\Tr T$ for $\Tr_{V_0}T$. If $\Tr T\ge 0$, then 
$$\kappa \lmax(T)\ge n_0 \lmax(T)\ge \Tr T.$$
The desired estimate follows. In the following, we may assume $\Tr T<0$. Let
$$q:=(p(x)-p(y))/2.$$
By Lemma~\ref{lem:H-formula},
\begin{equation}\label{eq:T-trace}
 \Tr T=-2(n_0+3)q.
\end{equation}
In particular, $q>0$.

Let
$$
 z={(x-y)}/{2} \quad \text{and}\quad v:={(x+y)}/{2}.
$$
Then $\ip{z}{v}=0$ and $|z|^2+|v|^2=1$.
By Theorem \ref{thm:spectral}, we may write
$$
 z=\sum_i\lambda_i c_i \quad \text{with}\quad
 \lambda_1\geq\cdots\geq\lambda_r,
 $$
for a Jordan frame $c_1,\ldots,c_r$. We use this Jordan frame throughout the remainder of the proof.

Expanding $p(v+z)-p(v-z)$ gives
$q=p(z)+3\ip{z}{v^2}$. Write
$s:=\lambda_1+\lambda_2$.
Applying both inequalities in Lemma~\ref{lem:main-algebra},
first to $z,z$ and then to the orthogonal pair $z,v$, yields
\begin{align*}
 q &\leq s|z|^2+\frac{3(r-1)}r s|v|^2 \\
 &= \frac{3(r-1)}{r} s -\frac{2r-3}{r} s |z|^2 \le \frac{3(r-1)}{r} s.
\end{align*}
In particular, $s>0$. The above also implies 
\begin{equation}\label{eq:q-upper}
3s-q\geq {(r-1)^{-1}}  {q} .
\end{equation}

By Lemma \ref{lem:H-formula}, we have for every unit $h\in V_0$
\begin{equation}\label{test-idenity}
 \ip{Th}{h}=12\ip{z}{h^2} -2q - 6 (\ip{x^2}{h}\ip{x}{h}-\ip{y^2}{h}\ip{y}{h} ) +3 (p(x){\ip{x}{h}}^2 - p(y){\ip{y}{h}}^2 ).
\end{equation}
In the following, we estimate $\lmax(T)$ by choosing specific unit test vectors $h$ in the above.

Since $V$ is simple, we know $\dim V_{12}=d>0$. Hence, we may
choose a unit $u\in V_{12}$. By Lemma~\ref{lem:rank-two}
and the properties of $V_{12}$,
$$
 u^2=\frac{c_1+c_2}{2},\quad
 z\circ u=\frac{s}{2}u,\quad\text{and} \quad 
 \ip{z}{u}=0.
$$
Since $x=v+z$ and $y=v-z$, it follows that
$$
 \ip{x}{u}=\ip{y}{u}=\ip{v}{u} \quad \text{and}\quad
 \ip{x^2-y^2}{u}=4\ip{z\circ v}{u}
 =2s\ip{v}{u}.
$$

Applying \eqref{test-idenity} with $h=u$ and 
using the above identities, we obtain
\begin{equation}\label{eq:first-test}
 \frac 12\lmax(T)\geq \frac12\ip{Tu}{u}
 =3s-q+3(q-2s)\ip{v}{u}^2.
\end{equation}
We divide the remainder of the proof into two cases.

\emph{Case 1.} $q\geq 2s$.
Then, by \eqref{eq:q-upper} and \eqref{eq:first-test}, we obtain
$$ 2^{-1} \lmax(T) \geq (r-1)^{-1}q \geq  (r+3)^{-1} q.$$
From this and \eqref{eq:T-trace}, the desired estimate is proved. 

\emph{Case 2.} $q < 2s$.
Recall that  $W_{12}:=\R(c_1-c_2)\oplus V_{12}$ is the traceless subspace of $E_1:=V_{11}\oplus V_{22}\oplus V_{12}$.
Since $\dim W_{12}=d+1\geq 2$, we may choose a unit
$w\in W_{12}$ perpendicular to $v$.
Since the components of $z$ and $z^2$ in $W_{12}$ are
$(\lambda_1-\lambda_2)(c_1-c_2)/2$ and
$(\lambda_1^2-\lambda_2^2)(c_1-c_2)/2$, respectively, we obtain
\begin{equation}\label{equ-5}
 \ip{z^2}{w}=s\ip{z}{w} \quad \text{and} \quad
 \sqrt2\,|\ip{z}{w}|\leq\lambda_1-\lambda_2.
\end{equation}
By $\tr z=0$, we have
$\lambda_1+(r-1)\lambda_2
\geq 0$. From this and the inequality in \eqref{equ-5}, we obtain
\begin{equation}\label{equ-6}
 \sqrt2\,|\ip{z}{w}|
 \leq {r}{(r-2)^{-1}}s\leq3s.
\end{equation}
Using the inequality in \eqref{equ-5} again, the fact $\lambda_1-\lambda_2\leq r\lambda_1/(r-1)$
and Lemma \ref{lem:spectral-bound}(i) yield 
\begin{equation}\label{equ-7}
     4\ip{z}{w}^2 
 \leq\frac{2r^2}{(r-1)^2}\lambda_1^2
 \leq\frac{2r}{r-1}|z|^2\leq3|z|^2.
\end{equation}

Applying Lemma~\ref{lem:rank-two} with $v\in w^\perp$
and the chosen $u\in V_{12}\subset E_1$ as above, we obtain
\begin{equation}\label{equ-8}
     2\sqrt2\,|\ip{v^2}{w}|
 \leq |v|^2-\ip{v}{u}^2,
\end{equation}
and also $w^2=(c_1+c_2)/2$.

By $\ip{v}{w}=0$, we have
$$\ip{x}{w}=\ip{z}{w}=-\ip{y}{w}.$$ 
Applying \eqref{test-idenity} with $h=w$ there and using the identity in \eqref{equ-5} and $x^2+y^2=2(z^2+v^2)$, we obtain
\begin{align*}
2^{-1}\lmax(T)\ge  2^{-1}\ip{Tw}{w}
 &=3s-q-6\ip{z}{w}\ip{v^2}{w}
             +3(q-2s)\ip{z}{w}^2\\
 &\geq3s-q-6|\ip{z}{w}|\,|\ip{v^2}{w}|
 -6s\ip{z}{w}^2.
 \end{align*}
Using \eqref{equ-6}, \eqref{equ-7}, and \eqref{equ-8}, we get
\begin{align*}
2^{-1}\lmax(T) &\geq 3s-q-\frac92s\bigl(|v|^2-\ip{v}{u}^2\bigr)
             -\frac92s|z|^2   \\
& = 3s-q-\frac92s\bigl( 1 -\ip{v}{u}^2\bigr).
\end{align*}

Multiply this by $2s-q>0$ and \eqref{eq:first-test} by $3s/2$, and then add them up.
The terms containing $\ip{v}{u}^2$ cancel, and we obtain
$$
 \left(\frac72s-q\right)\frac12\lmax(T)
 \geq q^2 -2qs +\frac{3 s^2}{2}.
$$
Dividing by the factor on the left hand side yields
$$
 \frac12\lmax(T)\geq\frac q6\geq\frac{q}{r+3}.
$$
Combining this with \eqref{eq:T-trace} proves the desired estimate.
\end{proof}

\section{Proof of Theorem \ref{thm:main}}
\label{sec:operator}

Fix the isometry $\iota:\R^{n_0}\to V_0$ in Theorem \ref{thm:main}.
For $u=U\circ\iota$, it holds that
\begin{equation}\label{eq:u-Hessian}
 D^2u(\xi)=\iota^*H(\iota\xi)\iota\quad\text{for }\xi\in \R^{n_0}\setminus\{0\},
\end{equation}
where
$\iota^*:V_0\to\R^{n_0}$ is the adjoint of $\iota$ defined by
$(\iota^*x)\cdot\xi=\ip{x}{\iota\xi}$ for $x\in V_0$ and $\xi\in\R^{n_0}$. Here $\cdot$ is the Euclidean inner product on $\R^{n_0}$. 
Clearly, $\iota^*=\iota^{-1}$.
The correspondence \eqref{eq:u-Hessian} allows us to pass the spectral properties of $H$ to those of $D^2 u$. In the rest of this section, we use
the Euclidean gradient $Du=(\partial_\alpha u)_\alpha$ and
Hessian matrix $D^2u=(\partial_{\alpha\beta}u)_{\alpha\beta}$.

The following lemma is a slight variant of
\cite[Lemma 3.1]{Chu-1}. For a symmetric matrix, $\tr$ denotes
its ordinary matrix trace.

\begin{lemma}\label{lem:extension}
Let $K\subset\Sym(n_0)$ be nonempty and compact.
Suppose that $\kappa\geq0$ and
\begin{equation}\label{eq:K-comparison}
 \tr(A-B)+\kappa\lmax(A-B)\geq 0 \quad \text{for all }A,B\in K.
\end{equation}
Then the operator $F:\Sym(n_0)\to\R$ defined by
\begin{equation}\label{eq:F-explicit}
 F(M):=\min_{A\in K}
 \bigl\{\tr(M-A)+\kappa\lmax(M-A)\bigr\},
 \quad M\in\Sym(n_0),
\end{equation}
satisfies $F=0$ on $K$ and \eqref{eq:ellipticity} with
$\lambda=1$ and $\Lambda=1+\kappa$.
\end{lemma}

\begin{proof}
The minimum exists by compactness of $K$.
For every $M\in\Sym(n_0)$, $N\geq 0$, and $A\in K$,
the Rayleigh quotient gives
$$
 0\leq\lmax(M+N-A)-\lmax(M-A)
 \leq \lmax(N)\leq\tr N.
$$
Hence,
\begin{align*}
 \tr N+\tr(M-A)+\kappa\lmax(M-A)
 &\leq\tr(M+N-A)+\kappa\lmax(M+N-A)\\
 &\leq(1+\kappa)\tr N +\tr(M-A)+\kappa\lmax(M-A).
\end{align*}
This implies \eqref{eq:ellipticity} with
$\lambda=1$ and $\Lambda=1+\kappa$.

Fix $B\in K$. Choosing $A=B$ in the definition gives
$F(B)\leq 0$. On the other hand, \eqref{eq:K-comparison}
gives $\tr(B-A)+\kappa\lmax(B-A)\geq0$ for every
$A\in K$, so $F(B)\geq 0$. Therefore $F(B)=0$. The lemma is proved. 
\end{proof}

Now we are ready to give the following:

\begin{proof}[Proof of Theorem~\ref{thm:main}]
The defining formula for $u=U\circ\iota$ shows that $u$ is
analytic away from zero and homogeneous of degree two. It follows that $u\in C^{1,1}(\R^{n_0})$.

Let $\kappa:=(r+3)(n_0+3)$ and
$$
 K:=\{D^2u(\xi)\mid |\xi|=1\}.
$$
The set $K$ is nonempty and compact. By \eqref{eq:u-Hessian},
Proposition~\ref{prop:separation} gives
$$
 \kappa\lmax\bigl(D^2u(\xi)-D^2u(\eta)\bigr)
 \geq\bigl|\tr\bigl(D^2u(\xi)-D^2u(\eta)\bigr)\bigr|
 \quad \text{for all }|\xi|=|\eta|=1.
$$
In particular, \eqref{eq:K-comparison} holds for every $A,B\in K$.
Apply Lemma~\ref{lem:extension} and define $F$ by
\eqref{eq:F-explicit}. Then $F=0$ on $K$ and $F$ satisfies
\eqref{eq:ellipticity} with the constants stated in the theorem.
From this and the homogeneity of $u$, we obtain
$$
 F(D^2u)=0\quad \text{in }\R^{n_0}\setminus\{0\}\quad \text{in the classical sense}.
$$

We next verify $F(D^2 u)=0$ across $0$ in the viscosity sense. We use the argument in
\cite[Proof of Theorem 1.1]{Chu}.
Let $\phi\in C^2$ touch $u$ from below at $0$.
Choose $\rho>0$ such that $u\geq\phi$ on $\overline{B_\rho}$
and $u(0)=\phi(0)$. Since $u\in C^1$, we have
$D\phi(0)=Du(0)=0$. For $\delta>0$ and
$0<\varepsilon<\delta\rho$, minimize
$$
 u(\xi)-\phi(\xi)+\delta|\xi|^2-\varepsilon\xi_1
 \quad\hbox{on }\overline{B_\rho}.
$$
This function is zero at $\xi=0$ and has gradient
$-\varepsilon e_1$ there. Its minimum is therefore negative.
It is bounded below by $\delta|\xi|^2-\varepsilon|\xi|$,
so any minimizer $\xi_\varepsilon$ satisfies
$$
 0<|\xi_\varepsilon|<{\varepsilon}/{\delta}<\rho.
$$
At this nonzero interior minimum, we have
$$
 D^2\phi(\xi_\varepsilon)-2\delta I
 \leq D^2u(\xi_\varepsilon).
$$
Ellipticity and the equation away from zero imply
$F(D^2\phi(\xi_\varepsilon)-2\delta I)\leq 0$.
Letting $\varepsilon\downarrow0$ and then $\delta\downarrow0$,
and using continuity of $F$, gives $F(D^2\phi(0))\leq 0$. This proves $F(D^2 u)\le 0$ near $0$ in the viscosity sense. A similar argument yields $F(D^2 u)\ge 0$ near $0$. Hence we have proved $F(D^2 u)=0$ in $\R^{n_0}$.

Finally, we check that $u$ is not twice differentiable at $0$.
Observe that
$u$ is odd and nonzero. To check the latter assertion,
choose a Jordan frame $c_1,\ldots,c_r$ and consider
$x_0=2c_1-c_2-c_3\in V_0$, which is possible because $r\ge 3$. By
Lemma \ref{lem:spectral-facts}(i),
$$
 U(x_0)=\frac{8-1-1}{\sqrt{4+1+1}}=\sqrt6 \ne 0.
$$
If $u$ were twice differentiable at zero, the homogeneity of $u$ and
$u(0)=Du(0)=0$ would force
$u(\xi)=2^{-1}\xi^{\mathsf T}D^2u(0)\xi$ for every $\xi$. The right side is even, contradicting
oddness of $u$ and $u\not\equiv 0$. This completes the proof of Theorem \ref{thm:main}.
\end{proof}

\end{document}